\RequirePackage{ifpdf}
\ifpdf 
\documentclass[pdftex]{sigma}
\else
\documentclass{sigma}
\fi

\begin{document}

\allowdisplaybreaks

\renewcommand{\PaperNumber}{016}

\FirstPageHeading

\ShortArticleName{On the Complex Symmetric and Skew-Symmetric Operators with a Simple Spectrum}

\ArticleName{On the Complex Symmetric and Skew-Symmetric\\ Operators with a Simple Spectrum}

\Author{Sergey M.~ZAGORODNYUK}

\AuthorNameForHeading{S.M.~Zagorodnyuk}

\Address{School of Mathematics and Mechanics, Karazin Kharkiv National University,\\
 4 Svobody Square, Kharkiv 61077, Ukraine}

\Email{\href{mailto:Sergey.M.Zagorodnyuk@univer.kharkov.ua}{Sergey.M.Zagorodnyuk@univer.kharkov.ua}}

\ArticleDates{Received December 14, 2010, in f\/inal form February 11, 2011;  Published online February 16, 2011}

\Abstract{In this paper we obtain necessary and suf\/f\/icient conditions for a linear bounded operator in
a Hilbert space $H$ to
have a three-diagonal complex symmetric  matrix with
non-zero elements on the f\/irst sub-diagonal
in an orthonormal basis in $H$. It is shown that a set of all such operators
is a proper subset of a set of all complex symmetric operators with a simple spectrum.
Similar necessary and suf\/f\/icient conditions are obtained
for a linear bounded operator in $H$ to
have a three-diagonal complex skew-symmetric  matrix with
non-zero elements on the f\/irst sub-diagonal
in an orthonormal basis in~$H$.}


\Classification{44A60}

\section{Introduction}
In last years an increasing interest is devoted to the subject of operators related to
bilinear forms in a Hilbert space (see~\cite{cit_1000_GP,cit_2000_GP2,cit_3000_Z} and
references therein), i.e.\ to the following forms:
\[ [x,y]_J := (x,Jy)_H,\qquad x,y\in H, \]
where $J$ is a conjugation and $(\cdot,\cdot)_H$ is the inner product in a Hilbert space $H$.
The conjugation~$J$ is an {\it antilinear} operator in $H$ such that $J^2 x = x$, $x\in H$,
and
\[ (Jx,Jy)_H = (y,x)_H,\qquad x,y\in H. \]
Recall that a linear operator $A$ in $H$ is said to be $J$-symmetric ($J$-skew-symmetric) if
\begin{gather}
\label{f1_2}
[Ax,y]_J = [x,Ay]_J,\qquad x,y\in D(A),
\end{gather}
or, respectively,
\begin{gather}
\label{f1_3}
[Ax,y]_J = -[x,Ay]_J,\qquad x,y\in D(A).
\end{gather}
If a linear bounded operator $A$ in a Hilbert space $H$ is $J$-symmetric ($J$-skew-symmetric)
for a conjugation $J$ in $H$, then $A$ is said to be complex symmetric (respectively complex skew-symmetric).
The matrices of complex symmetric (skew-symmetric) operators in certain bases of~$H$ are
complex symmetric (respectively skew-symmetric) semi-inf\/inite matrices.
Observe that for a bounded linear operator $A$ conditions~(\ref{f1_2}) and~(\ref{f1_3}) are
equivalent to conditions
\begin{gather}
\label{f1_4}
JAJ = A^*,
\end{gather}
and{\samepage
\[ JAJ = -A^*, \]
respectively.}

Recall that a bounded linear operator $A$ in a Hilbert space $H$ is said to have a simple spectrum if
there exists a vector $x_0\in H$ (cyclic vector) such that
\[ \overline{ \mathop{\rm Lin}\nolimits \{ A^k x_0,\ k\in\mathbb{Z}_+ \} } = H. \]
Observe that these operators are also called {\it cyclic operators}.

It is well known that a bounded self-adjoint operator with a simple spectrum
has a bounded semi-inf\/inite real symmetric three-diagonal (Jacobi) matrix in a certain orthonormal
basis (e.g.~\cite[Theorem 4.2.3]{cit_3500_Akh}).

The aim of our present investigation is to describe a class $C_+=C_+(H)$ ($C_-=C_-(H)$) of linear bounded operators in
a Hilbert space $H$, which
have three-diagonal complex symmetric (respectively skew-symmetric) matrices with
non-zero elements on the f\/irst sub-diagonal
in some orthonormal bases of $H$.
We obtain necessary and suf\/f\/icient conditions for a linear bounded operator in
a Hilbert space $H$ to belong to the class $C_{+}$ ($C_{-}$).
The class $C_{+}$ ($C_{-}$) is a subset of the class of all complex symmetric (respectively skew-symmetric)
operators in $H$ with a simple spectrum. Moreover, it is shown that $C_{+}(H)$ is a proper subset.

{\bf Notations.} As usual, we denote by $\mathbb{R}$, $\mathbb{C}$, $\mathbb{N}$, $\mathbb{Z}$, $\mathbb{Z}_+$
the sets of real numbers, complex numbers, positive integers, integers and non-negative integers,
respectively; $\mathop{\rm Im}\nolimits z = \frac{1}{2i}(z-\overline{z})$, $z\in \mathbb{C}$.
Everywhere in this paper, all Hilbert spaces are assumed to be separable. By~$(\cdot,\cdot)_H$ and~$\| \cdot \|_H$ we denote the scalar product and the norm in a Hilbert space $H$,
respectively. The indices may be omitted in obvious cases.
For a set $M$ in $H$, by $\overline{M}$ we mean the closure of $M$ in the norm $\| \cdot \|_H$. For
$\{ x_{n} \}_{n\in \mathbb{Z}_+}$, $x_{n}\in H$, we write
$\mathop{\rm Lin}\nolimits \{ x_{n} \}_{n\in \mathbb{Z}_+}$ for the set of linear combinations of elements
$\{ x_{n} \}_{n\in \mathbb{Z}_+}$.
The identity operator in $H$ is denoted by~$E_H$.
For an arbitrary linear operator $A$ in~$H$,
the operators $A^*$, $\overline{A}$, $A^{-1}$ mean its adjoint operator, its closure and its inverse
(if they exist). By $D(A)$ and $R(A)$ we mean the domain and the range of the operator $A$.
The norm of a bounded operator $A$ is denoted by $\| A \|$.
By $P^H_{H_1} = P_{H_1}$ we mean the operator of orthogonal projection in~$H$ on a subspace~$H_1$ in~$H$.

\section[The classes $C_\pm (H)$]{The classes $\boldsymbol{C_\pm (H)}$}

Let $\mathcal{M} = (m_{k,l})_{k,l=0}^\infty$, $m_{k,l}\in \mathbb{C}$, be a semi-inf\/inite complex matrix.
We shall say that $\mathcal{M}$ belongs to the class $\mathfrak{M}_3^+$, if and only if
the following conditions hold:
\begin{gather}
\label{f2_1}
m_{k,l} = 0,\qquad  k,l\in \mathbb{Z}_+,\qquad  |k-l|>1,
\\
\label{f2_2}
m_{k,l} = m_{l,k},\qquad k,l\in \mathbb{Z}_+,
\\
\label{f2_3}
m_{k,k+1} \not= 0,\qquad k\in \mathbb{Z}_+.
\end{gather}
We shall say that $\mathcal{M}$ belongs to the class $\mathfrak{M}_3^-$, if and only if
the conditions~(\ref{f2_1}), (\ref{f2_3}) hold and
\[
m_{k,l} = -m_{l,k},\qquad k,l\in \mathbb{Z}_+.\]
Let $A$ be a linear bounded operator in an inf\/inite-dimensional Hilbert space $H$.
We say that~$A$ belongs to the class $C_+ = C_+(H)$ ($C_- = C_-(H)$) if and only if there exists an
orthonormal basis~$\{ e_k \}_{k=0}^\infty$ in $H$ such that the matrix{\samepage
\begin{gather}
\label{f2_5}
\mathcal{M} = ( (Ae_l,e_k) )_{k,l=0}^\infty,
\end{gather}
belongs to $\mathfrak{M}_3^+$ (respectively to $\mathfrak{M}_3^-$).}

Let $y_0,y_1,\dots,y_n$ be arbitrary vectors in $H$, $n\in \mathbb{Z}_+$.
Set
\[
\Gamma(y_0,y_1,\dots,y_n) := \det ( (y_k,y_l)_H )_{k,l=0}^n. \]
Thus, $\Gamma(y_0,y_1,\dots,y_n)$ is the Gram determinant of vectors $y_0,y_1,\dots,y_n$.

The following theorem provides a description of the class $C_+(H)$.

\begin{theorem}
\label{t2_1}
Let $A$ be a linear bounded operator in an infinite-dimensional Hilbert space $H$.
The operator $A$ belongs to the class $C_+(H)$ if and only if the following conditions hold:
\begin{itemize}\itemsep=0pt
\item[$(i)$] $A$ is a complex symmetric operator with a simple spectrum;

\item[$(ii)$] there exists a cyclic vector $x_0$ of $A$ such that the following relations hold:
\begin{gather}
\label{f2_7}
\Gamma(x_0,x_1,\dots ,x_n,x_n^*) = 0,\qquad \forall\,  n\in \mathbb{N},
\end{gather}
where
\[
x_k = A^k x_0,\qquad x_k^* = \left(A^*\right)^k x_0,\qquad k\in \mathbb{N}; \]

and $Jx_0 = x_0$, for a conjugation $J$ in $H$ such that $JAJ=A^*$.
\end{itemize}
\end{theorem}

\begin{proof}
{\it Necessity.}
Let $H$ be an inf\/inite-dimensional Hilbert space and $A\in C_+(H)$.
Let $\{ e_k \}_{k=0}^\infty$ be an orthonormal basis in $H$ such that the matrix
$\mathcal{M} = (m_{k,l})_{k,l=0}^\infty$  belongs to
$\mathfrak{M}_3^+$, where $m_{k,l} = ( (Ae_l,e_k) )_{k,l=0}^\infty$.
Observe that
\begin{gather}
A e_0 = m_{0,0} e_0 + m_{1,0} e_1, \nonumber\\
A e_k = m_{k-1,k} e_{k-1} + m_{k,k} e_k + m_{k+1,k} e_{k+1},\qquad k\in \mathbb{N}.\label{f2_10}
\end{gather}
Suppose that
\[
e_r \in \mathop{\rm Lin}\nolimits \{ A^j e_0,\ 0\leq j\leq r \},\qquad 0\leq r\leq n, \]
for some $n\in \mathbb{N}$ (for $n=0$ it is trivial).
By~(\ref{f2_10}) we may write
\[ e_{n+1} = \frac{1}{ m_{n+1,n} } \left(
A e_n - m_{n-1,n} e_{n-1} - m_{n,n} e_n
\right) \in \mathop{\rm Lin}\nolimits \{ A^j e_0,\ 0\leq j\leq n+1 \}. \]
Here $m_{-1,0}:=0$ and $e_{-1} := 0$.
By induction we conclude  that
\begin{gather}
\label{f2_12}
e_r \in \mathop{\rm Lin}\nolimits \{ A^j e_0,\ 0\leq j\leq r \},\qquad r\in \mathbb{Z}_+.
\end{gather}
Therefore
$\overline{ \mathop{\rm Lin}\nolimits \{ A^j e_0,\ j\in \mathbb{Z}_+ \} } = H$, i.e.\ the operator
$A$ has a simple spectrum and $e_0$ is a cyclic vector of~$A$.

Consider the following conjugation:
\[
J \sum_{k=0}^\infty x_k e_k = \sum_{k=0}^\infty \overline{x_k} e_k,\qquad x=\sum_{k=0}^\infty x_k e_k\in H. \]

Observe that
\[ [Ae_k,e_l]_J = (Ae_k,e_l) = m_{l,k} = m_{k,l} = (Ae_l,e_k) = [Ae_l,e_k]_J,\qquad k,l\in \mathbb{Z}_+. \]
By linearity of the $J$-form $[\cdot,\cdot]_J$ in the both arguments we get
\[ [Ax,y]_J = [Ay,x]_J,\qquad x,y\in H. \]
Thus, the operator $A$ is $J$-symmetric and relation~(\ref{f1_4}) holds. Notice that
$Je_0 = e_0$. It remains to check if relation~(\ref{f2_7}) holds.
Set
\[
H_r := \mathop{\rm Lin}\nolimits \{ A^j e_0,\ 0\leq j\leq r \},\qquad r\in \mathbb{Z}_+.
\]
By~(\ref{f2_12}) we see that $e_0,e_1,\dots ,e_r\in H_r$, and therefore
$\{ e_j \}_{j=0}^r$ is an orthonormal basis in $H_r$ ($r\in \mathbb{Z}_+$).
Since $Je_j = e_j$, $j\in \mathbb{Z}_+$, we have
\[
JH_r \subseteq H_r,\qquad r\in \mathbb{Z}_+.
\]
Then
\[ \left( A^* \right)^r e_0 = \left( JAJ \right)^r e_0 = J A^r Je_0 = J A^r e_0\in H_r,\qquad r\in \mathbb{Z}_+. \]
Therefore vectors
$e_0, Ae_0,\dots ,A^r e_0, \left( A^* \right)^r e_0$,
are linearly dependent and their Gram determinant is equal to zero. Thus, relation~(\ref{f2_7}) holds
with $x_0 = e_0$.

{\it Sufficiency.}
Let $A$ be a bounded operator in a Hilbert space $H$ satisfying conditions~$(i)$, $(ii)$ in the statement
of the theorem.
For the cyclic vector $x_0$ we set
\[
H_r := \mathop{\rm Lin}\nolimits \{ A^j x_0,\ 0\leq j\leq r \},\qquad r\in \mathbb{Z}_+. \]

Observe that
\begin{gather}
\label{f2_17}
A^{r+1} x_0 \notin H_r,\qquad r\in \mathbb{Z}_+.
\end{gather}
In fact, suppose that for some $k\in \mathbb{N}$, we have
\[
A^{r+j} x_0 \in H_r,\qquad 1\leq j\leq k.
\]
Then
\begin{gather*}
A^{r+k+1} x_0 = A A^{r+k}x_0 = A \sum_{t=0}^r \alpha_{r,k;t} A^t x_0   = \sum_{t=0}^r \alpha_{r,k;t} A^{t+1} x_0\in H_r,\qquad \alpha_{r,k;t}\in \mathbb{C}.\nonumber
\end{gather*}
By induction we obtain
\[
A^{r+j} x_0 \in H_r,\qquad j\in \mathbb{Z}_+.\]
Therefore $H=H_r$. We obtain a contradiction since $H$ is inf\/inite-dimensional.

Let us apply the Gram--Schmidt orthogonalization method to the sequence
$x_0, Ax_0, A^2 x_0,\dots $. Namely, we set
\[
g_0 = \frac{x_0}{ \| x_0 \|_H },\qquad g_{r+1} = \frac{ A^{r+1}x_0 - \sum\limits_{j=0}^r (A^{r+1}x_0,g_j)_H g_j }
{ \Big\| A^{r+1}x_0 - \sum\limits_{j=0}^r (A^{r+1}x_0,g_j)_H g_j \Big\|_H },\qquad r\in \mathbb{Z}_+.
\]
By construction we have
\[
H_r = \mathop{\rm Lin}\nolimits \{ g_j,\ 0\leq j\leq r \},\qquad r\in \mathbb{Z}_+. \]
Therefore $\{ g_j \}_{j=0}^r$ is an orthonormal basis in $H_r$ ($r\in \mathbb{Z}_+$) and
$\{ g_j \}_{j\in \mathbb{Z}_+}$ is an orthonormal basis in $H$.

From~(\ref{f2_7}) and~(\ref{f2_17}) we conclude that
\[
J A^n x_0 = J A^n J x_0 = \left( A^* \right)^n x_0 \in H_n,\qquad n\in \mathbb{Z}_+. \]
Therefore
\begin{gather}
\label{f2_23}
J H_r \subseteq H_r,\qquad r\in \mathbb{Z}_+.
\end{gather}
Let
\[ J g_r = \sum_{j=0}^r \beta_{r,j} g_j,\qquad \beta_{r,j}\in \mathbb{C},\qquad r\in \mathbb{Z}_+. \]
Using properties of the conjugation and relation~(\ref{f2_23}) we get
\[ \beta_{r,j} = ( J g_r, g_j )_H = ( J g_r, JJ g_j )_H = \overline{ ( g_r, J g_j )_H } = 0, \]
for $0\leq j\leq r-1$.
Therefore
\[
J g_r = \beta_{r,r} g_r,\qquad \beta_{r,r}\in \mathbb{C},\qquad r\in \mathbb{Z}_+.
\]
Since $\| g_r \|^2 = \| Jg_r \|^2 = | \beta_{r,r} |^2 \| g_r \|^2$, we
have
\[
\beta_{r,r} = e^{i\varphi_r},\qquad \varphi_r\in [0,2\pi),\qquad  r\in \mathbb{Z}_+.
\]
Set
\[
e_r :=  e^{i \frac{\varphi_r}{2} } g_r,\qquad r\in \mathbb{Z}_+.
\]
Then $\{ e_j \}_{j=0}^r$ is an orthonormal basis in $H_r$ ($r\in \mathbb{Z}_+$) and
$\{ e_j \}_{j\in \mathbb{Z}_+}$ is an orthonormal basis in~$H$.
Observe that
\[
J e_r =  J e^{i \frac{\varphi_r}{2} } g_r = e^{-i \frac{\varphi_r}{2} } J g_r =
e^{i \frac{\varphi_r}{2} } g_r = e_r,\qquad r\in \mathbb{Z}_+.
\]
Def\/ine the matrix $\mathcal{M} = (m_{k,l})_{k,l=0}^\infty$ by~(\ref{f2_5}). Notice that
\[ m_{k,l} = (Ae_l,e_k)_H = [Ae_l,e_k]_J = [e_l,Ae_k]_J = [Ae_k,e_l]_J = (Ae_k,e_l)_H = m_{l,k}, \]
where $k,l\in \mathbb{Z}_+$, and therefore $\mathcal{M}$ is complex symmetric.

If $l\geq k+2$ ($k,l\in \mathbb{Z}_+$), then
\[ m_{k,l} = (Ae_l,e_k)_H = [Ae_l,e_k]_J = [e_l,Ae_k]_J = (e_l, JAe_k)_H = 0, \]
since $JAe_k \in H_{k+1}\subseteq H_{l-1}$, and $e_l \in H_l\ominus H_{l-1}$.
Therefore $\mathcal{M}$ is three-diagonal.

Since $e_r\in H_r$, using the def\/inition of $H_r$ we get
\[
A e_r\subseteq H_{r+1},\qquad r\in \mathbb{Z}_+. \]
Observe that
\[
A e_r \notin H_r,\qquad r\in \mathbb{Z}_+. \]
In fact, in the opposite case we get
\[ A e_j \in H_r,\qquad 0\leq j\leq r, \]
and $AH_r \subseteq H_r$. Then $A^k x_0\in H_r$, $k\in \mathbb{Z}_+$, and
$H=H_r$. This is a contradiction since $H$ is an inf\/inite-dimensional space.

Hence, we may write
\[
A e_r = \sum_{j=0}^{r+1} \gamma_{r,j} e_j,\qquad \gamma_{r,j}\in \mathbb{C},\qquad \gamma_{r,r+1}\not=0. \]
Observe that
\[ m_{r+1,r} = (Ae_r,e_{r+1})_H = \gamma_{r,r+1}\not= 0,\qquad r\in \mathbb{Z}_+. \]
Thus, $\mathcal{M}\in \mathfrak{M}_3^+$ and $A\in C_+(H)$.
\end{proof}

\begin{remark}
Condition~$(ii)$ of the last theorem may be replaced by the following condition which
does not use a conjugation $J$:
\begin{itemize}\itemsep=0pt
\item[$(ii)^*$] there exists a cyclic vector $x_0$ of $A$ such that the following relations hold:
\begin{gather}
\label{f2_31}
\Gamma(x_0,x_1,\dots ,x_n,x_n^*) = 0,\qquad \forall\, n\in \mathbb{N},
\end{gather}
where
\[
x_k = A^k x_0,\qquad x_k^* = \left(A^*\right)^k x_0,\qquad k\in \mathbb{N}, \]
and the following operator:
\begin{gather}
\label{f2_33}
L \sum_{k=0}^\infty \alpha_k A^k x_0 := \sum_{k=0}^\infty \overline{\alpha_k} \left( A^* \right)^k x_0,\qquad
\alpha_k\in \mathbb{C},
\end{gather}
where all but f\/inite number of coef\/f\/icients $\alpha_k$ are zeros,
is a bounded operator in $H$ which extends by continuity to a conjugation in $H$.
\end{itemize}

Let us show that conditions $(i)$, $(ii)$ $\Leftrightarrow$ conditions $(i)$, $(ii)^*$.

The necessity is obvious since the conjugation $J$ satisf\/ies relation~(\ref{f2_33})
(with $J$ instead of $L$).

{\it Sufficiency.}
Let conditions $(i)$, $(ii)^*$ be satisf\/ied. Notice that
\begin{gather}
LA A^k x_0 = L A^{k+1} x_0 = \left( A^* \right)^{k+1} x_0, \nonumber\\
A^* L A^k x_0 = A^* \left( A^* \right)^k x_0 = \left( A^* \right)^{k+1} x_0,\qquad k\in \mathbb{Z}_+. \nonumber
\end{gather}
By continuity we get $ LA = A^* L$.
Then condition~$(ii)$ holds with the conjugation~$L$.
\end{remark}

\begin{remark}
Notice that conditions~(\ref{f2_31}) may be written in terms of
the coordinates of $x_0$ in an arbitrary orthonormal basis $\{ u_n \}_{n=0}^\infty$ in $H$:
\[
x_0 = \sum_{n=0}^\infty x_{0,n} u_n,\qquad
A^k x_0 = \sum_{n=0}^\infty x_{0,n} A^k u_n,\qquad
\left( A^* \right)^k x_0 = \sum_{n=0}^\infty x_{0,n} \left( A^* \right)^k u_n. \]
By substitution these equalities in relation~(\ref{f2_7}) we get some algebraic
equations with respect to the coordinates $x_{0,n}$.
If cyclic vectors of $A$ are unknown, one can use numerical methods
to f\/ind approximate solutions of these equations.
Then there should be cyclic vectors of $A$ among these solutions.
\end{remark}

The following theorem gives an analogous description for the class $C_-(H)$.
\begin{theorem}
\label{t2_2}
Let $A$ be a linear bounded operator in an infinite-dimensional Hilbert space $H$.
The operator $A$ belongs to the class $C_-(H)$ if and only if the following conditions hold:
\begin{itemize}\itemsep=0pt
\item[$(i)$] $A$ is a complex skew-symmetric operator with a simple spectrum;

\item[$(ii)$] there exists a cyclic vector $x_0$ of $A$ such that the following relations hold:
\begin{gather}
\label{f2_36}
\Gamma(x_0,x_1,\dots ,x_n,x_n^*) = 0,\qquad \forall \, n\in \mathbb{N},
\end{gather}
where
\begin{gather}
\label{f2_37}
x_k = A^k x_0,\qquad x_k^* = \left(A^*\right)^k x_0,\qquad k\in \mathbb{N};
\end{gather}
and $Jx_0 = x_0$, for a conjugation $J$ in $H$ such that $JAJ=-A^*$.
\end{itemize}
\end{theorem}

Condition~$(ii)$ of this theorem may be replaced by the following condition:
\begin{itemize}\itemsep=0pt
\item[$(ii)^*$] there exists a cyclic vector $x_0$ of $A$ such that relations~(\ref{f2_36}), (\ref{f2_37}) hold
and the following operator:
\[
L \sum_{k=0}^\infty \alpha_k A^k x_0 := \sum_{k=0}^\infty (-1)^k \overline{\alpha_k} \left( A^* \right)^k x_0,\qquad
\alpha_k\in \mathbb{C}, \]
where all but f\/inite number of coef\/f\/icients $\alpha_k$ are zeros,
is a bounded operator in $H$ which extends by continuity to a conjugation in $H$.
\end{itemize}
The proof of the latter facts is similar and essentially the same as for the case of $C_+(H)$.

The following example shows that condition~$(ii)$ (or~$(ii)^*$) can not be removed from Theorem~\ref{t2_1}.

\begin{example}
\label{e2_1}
Let $\sigma(\theta)$ be a non-decreasing left-continuous bounded function on $[0,2\pi]$
with an inf\/inite number of points of increase and
such that
\begin{gather}
\label{f2_39}
\int_0^{2\pi} \ln \sigma' (\theta) d\theta = -\infty.
\end{gather}
Consider the Hilbert space $L^2([0,2\pi],d\sigma)$ of (classes of equivalence of) complex-valued func\-tions~$f(\theta)$
on $[0,2\pi]$ such that
\[
\| f \|^2_{L^2([0,2\pi],d\sigma)} :=
\left(
\int_0^{2\pi} |f(\theta)|^2 d\sigma(\theta)
\right)^{ \frac{1}{2} }
<\infty. \]

The condition~(\ref{f2_39}) provides that algebraic polynomials of $e^{i\theta}$ are dense in
$L^2([0,2\pi],d\sigma)$ \mbox{\cite[p.~19]{cit_4000_G}}.
Therefore the operator
\[
U f(\theta) = e^{i\theta} f(\theta),\qquad f\in L^2([0,2\pi],d\sigma), \]

is a cyclic unitary operator in $H$, with a cyclic vector $f_0(\theta)=1$.
Set
\[
J f(\theta) = \overline{ f(\theta) },\qquad f\in L^2([0,2\pi],d\sigma). \]

Then
\[ JUJ f(\theta) = J e^{i\theta} \overline{ f(\theta) } = e^{-i\theta} f(\theta) = U^{-1} f(\theta) =
U^* f(\theta). \]
Thus, $U$ is a complex symmetric operator with a simple spectrum and condition~$(i)$ of
Theorem~\ref{t2_1} is satisf\/ied.

However, $U\notin C_+(H)$. In fact, suppose to the contrary that there exists an orthonormal
basis $\{ e_j \}_{j\in \mathbb{Z}_+}$ such that the corresponding matrix $\mathcal{M} =
(m_{k,l})_{k,l=0}^\infty$ from~(\ref{f2_5})
belongs to the class~$\mathfrak{M}_3^+$.
Since $U$ is unitary, we have
\[ \mathcal{E} = \mathcal{M} \mathcal{M}^*, \]
with the usual rules of matrix operations, $\mathcal{E} = (\delta_{k,l})_{k,l=0}^\infty$.
However, the direct calculation shows that the element of the matrix $\mathcal{M} \mathcal{M}^*$
in row~0, column~2 is equal to $m_{0,1} \overline{m_{2,1}} \not= 0$.
We obtained a contradiction. Thus, $U\notin C_+(H)$.
Consequently, {\it condition~$(ii)$ in Theorem~{\rm \ref{t2_1}} is essential and can not be removed.}
\end{example}

\begin{proposition}
\label{p2_1}
Let $H$ be an arbitrary infinite-dimensional Hilbert space. The class
$C_+(H)$ is a proper subset of the set of all complex symmetric
operators with a simple spectrum in $H$.
\end{proposition}
\begin{proof}
Consider an arbitrary inf\/inite-dimensional Hilbert space $H$.
Let $V$ be an arbitrary unitary operator which maps $L^2([0,2\pi],d\sigma)$ (see Example~\ref{e2_1})
onto $H$. Then $\widehat U := VUV^{-1}$ is a unitary operator in $H$ with a simple spectrum
and it has a cyclic vector $\widehat x_0 := V 1$.
Since $JUJ = U^*$, we get
\begin{gather*}
JV^{-1} \widehat U VJ = V^{-1} \widehat U^* V,\qquad
VJV^{-1} \widehat U VJV^{-1} = \widehat U^*.
\end{gather*}
Observe that $\widehat J := VJV^{-1}$ is a conjugation in $H$. Therefore $\widehat U$ is
a complex symmetric operator in $H$.
Suppose that $\widehat U\in C_+(H)$.
Let $\mathcal{F} = \{ f_k \}_{k=0}^\infty$ be
an orthonormal basis in $H$ such that the matrix $M = (m_{k,l})_{k,l=0}^\infty$,
$m_{k,l} = (\widehat U f_l,f_k)_H$, belongs to $\mathfrak{M}_3^+$.
Observe that $\mathcal{G} = \{ g_k \}_{k=0}^\infty$, $g_k := V^{-1} f_k$, is an orthonormal
basis in $L^2([0,2\pi],d\sigma)$ and
\begin{gather*}
(U g_l,g_k)_{L^2([0,2\pi],d\sigma)} = (V^{-1} \widehat U V g_l,g_k)_{L^2([0,2\pi],d\sigma)} =
(\widehat U f_l,f_k)_H   = m_{k,l},\qquad k,l\in \mathbb{Z}_+.
\end{gather*}
Therefore $U\in C_+(L^2([0,2\pi],d\sigma))$. This is a contradiction with Example~\ref{e2_1}.
Consequently, we have $\widehat U\notin C_+(H)$.

On the other hand, the class $C_+(H)$ is non-empty, since an arbitrary matrix from $\mathcal{M}_3^+$
with bounded elements def\/ine an
operator $B$ in $H$ which have this matrix in an arbitrary f\/ixed orthonormal basis in $H$.
\end{proof}

\begin{remark}
The classical Jacobi matrices
are closely related to orthogonal polynomials~\cite{cit_3500_Akh}.
Let us indicate some similar relations for the class $\mathfrak{M}_3^+$.
Choose an arbitrary
$\mathcal{M} = (m_{k,l})_{k,l=0}^\infty\in \mathfrak{M}_3^+$, where $m_{k,l}\in \mathbb{C}$.
Let $\{ p_n(\lambda) \}_{n=0}^\infty$, $\deg p_n=n$, $p_0(\lambda)=1$, be a sequence of polynomials def\/ined
recursively by the following relation:
\begin{gather}
\label{f2_43}
m_{n,n-1} p_{n-1}(\lambda) + m_{n,n} p_n(\lambda) + m_{n,n+1} p_{n+1}(\lambda) = \lambda p_n(\lambda),\qquad
n=0,1,2,\dots ,
\end{gather}
where $m_{0,-1}:=1$, $p_{-1} := 0$.
Set $c_n = m_{n,n+1}$, $b_n = m_{n,n}$, $n\in \mathbb{Z}_+$; and $c_{-1}:=1$.
By~(\ref{f2_2}), (\ref{f2_43}) we get
\begin{gather}
\label{f2_44}
c_{n-1} p_{n-1}(\lambda) + b_n p_n(\lambda) + c_n p_{n+1}(\lambda) = \lambda p_n(\lambda),\qquad
n=0,1,2,\dots .
\end{gather}
Let $p_n(\lambda) = \mu_n \lambda^n + \cdots $, $\mu_n\in \mathbb{C}$, $n\in \mathbb{Z}_+$.
Comparing coef\/f\/icients by $\lambda^{n+1}$ in~(\ref{f2_44}) we get
\[ \mu_{n+1} = \frac{1}{c_n} \mu_n,\qquad n\in \mathbb{Z}_+. \]
By induction we see that
\[
\mu_n = \left( \prod_{j=0}^{n-1} c_j \right)^{-1},\qquad n\in \mathbb{N},\qquad \mu_0=1. \]

Set
\[
P_n(\lambda) = \prod_{j=0}^{n-1} c_j p_n(\lambda),\qquad n\in \mathbb{N},\qquad P_0(\lambda)=1,\qquad P_{-1}(\lambda)=0. \]

Multiplying the both sides of~(\ref{f2_44}) by $\prod\limits_{j=0}^{n-1} c_j$, $n\geq 1$, we obtain:
\[
c_{n-1}^2 P_{n-1}(\lambda) + b_n P_n(\lambda) + P_{n+1}(\lambda) = \lambda P_n(\lambda),\qquad
n=0,1,2,\dots . \]

By~Theorem~6.4 in~\cite{cit_5000_Ch} there exists a complex-valued function $\phi$ of
bounded variation on $\mathbb{R}$ such that
\[
\int_\mathbb{R} P_m(\lambda) P_n(\lambda) d\phi(\lambda) = \left( \prod_{j=0}^{n-1} c_j \right)^2 \delta_{m,n},\qquad
m,n\in \mathbb{Z}_+. \]
Therefore we get
\[
\int_\mathbb{R} p_m(\lambda) p_n(\lambda) d\phi(\lambda) = \delta_{m,n},\qquad m,n\in \mathbb{Z}_+. \]

Polynomials $\{ p_n(\lambda) \}_{n=0}^\infty$ were used in~\cite{cit_5500_G,cit_6000_Z} to state and solve
the direct and inverse spectral problems for matrices from $\mathfrak{M}_3^+$.
Analogs of some facts of the Weyl discs theory were obtained for the case of matrices from~$\mathfrak{M}_3^+$
with additional assumptions~\cite{cit_7000_Z}:
$m_{n,n+1} > 0$, $n\in \mathbb{Z}_+$, and
\[
m_{n,n}\in \mathbb{C}: \ \ r_0 \leq \mathop{\rm Im}\nolimits m_{n,n} \leq r_1, \]
for some $r_0,r_1\in \mathbb{R}$, $n\in \mathbb{Z}_+$.

On the other hand, the direct and inverse spectral problems for matrices from $\mathfrak{M}_3^-$ were
investigated in~\cite{cit_8000_Z}.

Probably, some progress in the spectral theory of complex symmetric and skew-symmetric operators
would provide some additional information about corresponding polynomials and vice versa.
\end{remark}

\subsection*{Acknowledgements}

The author is grateful to referees for their comments and suggestions.

\pdfbookmark[1]{References}{ref}
\LastPageEnding

\end{document}